\documentclass[12pt,a4paper]{article}
\usepackage{amstext, amssymb, amscd, amsmath, amsthm, tabls}
\usepackage[all]{xy}
\usepackage{epsfig}

\newcommand{\ra}{{\rightarrow}}
\newcommand{\lra}{{\longrightarrow}}

\newtheorem{thm}{Theorem}[section]
\newtheorem{prop}[thm]{Proposition}

\newtheorem{lemma}[thm]{Lemma}

\theoremstyle{definition}
\newtheorem{defin}[thm]{Definition}

\theoremstyle{remark}
\newtheorem*{rem}{Remark}

\theoremstyle{remark}

\title{Symplectic geography in dimension $8$.}

\author{Federica Pasquotto}

\begin{document}

\maketitle

\medskip

\noindent
{\bf Abstract.}
We show that in dimension $8$ the geography of symplectic manifolds 
does not differ from that of almost complex ones.

\begin{section}{Introduction}
The Riemann-Roch theorem implies necessary 
and sufficient conditions for a given system of integer numbers to appear
as the system of Chern numbers of an almost complex manifold.
In dimension $8$ these relations are (see \cite{hi60})
\begin{eqnarray}\label{congrel}
-c_4+c_1c_3+3c_2^2+4c_2c_1^2-c_1^4 & \equiv & 0\ (\textrm{mod}\ 720)
\nonumber\\
2c_1^4+c_1^2c_2 & \equiv & 0\ (\textrm{mod}\ 12)\\
-2c_4+c_1c_3 & \equiv & 0\ (\textrm{mod}\ 4).\nonumber
\end{eqnarray}
A quintuple of integer numbers satisfying the above congruence relations
is called admissible.
We ask ourselves which admissible quintuples may appear as the
system of Chern numbers of a closed, connected, symplectic $8$-dimensional 
manifold. The problem of determining which systems of integer numbers
admit a symplectic realisation is known as symplectic geography. 

We prove: 
\begin{thm}\label{main}
Any ordered quintuple of integers which arises as the system of Chern 
numbers of an almost complex $8$-dimensional manifold can also be realised 
by a closed, connected, symplectic $8$-manifold.
\end{thm}  
The analogous result was proved by Halic \cite{ha99} in dimension $6$.
His method relies on symplectic constructions such as blow-up and connected 
symplectic sum. We follow his strategy of proof and make use of a general 
formula for computing the Chern classes of blow-up, which was known for 
algebraic varieties \cite{lasc75}
and whose proof can easily be modified so that it applies
to symplectic manifolds. Moreover, we apply Donaldson's theorem about 
existence of symplectic submanifolds to the total space of some symplectic 
fibrations.

\end{section}

\begin{section}{The eight-dimensional case.} 
\begin{subsection}{Congruence relations in dimension eight.}

Suppose we are given a
quintuple of integer numbers $(c_4, c_1c_3, c_2^2, c_1^2c_2, c_1^4)$
satisfying the system of congruence relations (\ref{congrel}).
Then there exist integers $(a, j, k, m, b)$ such that 
\begin{eqnarray}\label{congpar}
a & = & c_4 \nonumber \\
720j & =& -c_4+c_1c_3+3c_2^2+4c_2c_1^2-c_1^4 \nonumber \\
12k & = & 2c_1^4+c_1^2c_2 \\
4m & = & -2c_4+c_1c_3 \nonumber \\
b & = & c_1^4 \nonumber
\end{eqnarray}
and the above system is equivalent to 
\begin{eqnarray*}
c_4 & = & a\\
c_1c_3 & = & 4m+2a\\
c_1^4 & = & b\\
c_1^2c_2 & = & 12k-2b\\
3c_2^2 & = & 720j-a-4m-48k+9b.
\end{eqnarray*}
From this we see that there is a one-to-one correspondence between quintuples 
of integers satisfying (\ref{congrel}) and quintuples $(a, b, j, k, m)$ 
subject to the condition $a+m\equiv 0\ (\textrm{mod}\ 3)$. 

The result which we want to prove can be summarised by saying that all
 quintuples $(a, b, j, k, m)$ satisfying the condition 
$a+m\equiv 0\ (\textrm{mod}\ 3)$ admit a symplectic realisation.
 
\begin{thm}\label{geopar}
Given a quintuple of integer numbers $(a, j, k, m, b)$, subject to the 
additional condition $a+m\equiv 0\ (\textrm{mod}\ 3)$, there exists a 
closed, connected, symplectic $8$-dimensional  manifold $M$ such that 
the given parameters are related to the Chern numbers of $M$ by the system 
of equations {\rm \eqref{congpar}}.
\end{thm}

In view of the above correspondence, the proposition immediately implies
that the congruence relations (\ref{congrel}) are not only necessary, but
also sufficient for a given quintuple of integer numbers to occur as the
system of Chern numbers of a closed connected $8$-dimensional symplectic 
manifold.
This is equivalent to 
the statement of Theorem \ref{main}.
\end{subsection}

\begin{subsection}{Behaviour of the parameters under blow-up.}

Let $\widetilde{M}$ denote the blow-up of the symplectic manifold 
$M$ along a symplectic submanifold $N$. In \cite{lasc75} the authors
prove a formula for the Chern classes of the blow-up of
an algebraic variety along a subvariety. This formula in fact holds
for symplectic manifolds as well and gives an expression for the total 
Chern class of $\widetilde{M}$ in terms of those of $M$ and $N$.
From the corresponding expressions for the individual Chern classes,
a straightforward computation shows that 
in dimension $8$, the parameters $(a, m, j, k, b)$ transform under blow-up
as follows \cite{pa04}.
\begin{itemize}
\item Blow-up at a point:
\begin{eqnarray}\label{parpoint}
a' & = & a+3\nonumber\\
4m' & = & 4m\nonumber\\
720j' & = & 720j\\
12k' & = & 12k-180\nonumber\\
b' & = & b-81\nonumber
\end{eqnarray}
\item Blow-up along a symplectically embedded curve $C$ of genus $g$ and 
with normal bundle $\nu C$: 
\begin{eqnarray}\label{parcurve}
a' & = & a+4(1-g)\nonumber\\
4m' & = & 4m-4(1-g)\nonumber\\
720j' & = & 720j\\
12k' & = & 12k-144(1-g)-36\langle c_1(\nu_C), [C]\rangle\nonumber \\
b' & = & b-64(1-g)-16\langle c_1(\nu_C), [C]\rangle \nonumber
\end{eqnarray}
\item Blow up along a symplectic four-dimensional submanifold $X$, with  
normal bundle $\nu X$:
\begin{eqnarray}\label{parfourdim}
a' & = & a+c_2[X]\nonumber\\
4m' & = & 4m+c_1^2[X]-3c_2[X]\nonumber\\
720j' & = & 720j\\
12k' & = & 12k-13c_1^2[X]-c_2[X]-18\langle c_1(X)c_1(\nu_X), [X]\rangle \nonumber\\
     & & -6\langle c_1^2(\nu_X), [X]\rangle \nonumber\\
b' & = & b-6c_1^2[X]-8\langle c_1(X)c_1(\nu_X), [X]\rangle 
-3\langle c_1^2(\nu_X),
[X]\rangle \nonumber \\
  & & +\langle c_2(\nu_X), [X]\rangle \nonumber
\end{eqnarray}
\end{itemize}
Notice that the 
parameter $j$ defined in (\ref{congpar}) is invariant under blow up. 

\end{subsection}
\end{section}

\medskip

\begin{section}{Building blocks.}

\begin{subsection}{Elliptic surfaces.}

For more details about the content of this paragraph 
we refer to \cite{gost99}.     
Let $E(1)$ denote the manifold 
$\mathbb{C}\mathbb{P}^2\#\,
9\, \overline{\mathbb{C}\mathbb{P}^2}$, equipped with an elliptic
fibration. 
The homology of $E(1)$ is the homology of a connected sum and it is given by 
$$H_2(E(1))\cong H_2(\mathbb{C}\mathbb{P}^2\#\, 9\,
\overline{\mathbb{C}\mathbb{P}^2})\cong H_2(\mathbb{C}\mathbb{P}^2;
\mathbb{Z})\oplus 9 H_2(\overline{\mathbb{C}\mathbb{P}^2}; \mathbb{Z}).$$ 
Let $h$ be the positive generator of  
$H_2(\mathbb{C}\mathbb{P}^2; \mathbb{Z})$,
$e_i$ the exceptional sphere of the i-th blow up, $i=1, \ldots , 9$, 
$f=3h-\sum_{i=1}^9 e_i$ the class of a regular fibre:
a basis for $H_2(E(1))$ is given by 
$\langle f, e_9, e_1-e_2, \ldots , e_7-e_8, -h+e_6+e_7+e_8\rangle $
and the intersection matrix with respect to this basis is 
\begin{eqnarray*}
\left [ 
\begin{array}{cc}
0 & 1\\
1 & -1
\end{array}
\right ]
\oplus (-E_8).
\end{eqnarray*}
The first Chern class of $E(1)$ is $PD(3h-\sum e_i)=PD(f)$, 
which gives for the Chern numbers 
the values $c_1^2[E(1)]=0$ and $c_2[E(1)]=12$.  


Inductively, one can perform symplectic sums along regular fibres and
define $E(n+1):=E(n)\#_F E(1)=\#_F^{n+1} E(1)$. 
This admits a basis
for the second homology group with corresponding intersection matrix 
\begin{eqnarray}\label{Enint}
n(-E_8)\oplus
2(n-1)\left [ 
\begin{array}{cc}
0 & 1\\
1 & -2
\end{array}
\right ]\oplus
\left [
\begin{array}{cc}
0 & 1\\
1 & -n
\end{array}
\right ]
.
\end{eqnarray}
We denote the elements of this basis by
\begin{equation*}
\langle\tau_{ij}, i=1,\ldots , n,\, j=1,\ldots ,8;\, \alpha_k, \beta_k, k=1, 
\ldots ,2(n-1);\, f, \sigma \rangle ,
\end{equation*} 
where $\sigma$ denotes the class of a section of $E(n)$, 
which is obtained by pasting together $n$ sections of $E(1)$.

The first Chern class is $c_1(E(n))=(2-n)f$, so the Chern numbers are given by 
\begin{eqnarray*}
\left\{
\begin{array}{l}
c_1^2[E(n)]=0\\
c_2[E(n)]=12n.
\end{array}
\right.
\end{eqnarray*}

\begin{defin}\label{nucleus}
The \underline{nucleus} $N(n)$ of the elliptic surface $E(n)$
consists of a neighbourhood of the union of a singular 
fibre and a section of the fibration.
\end{defin}

If we consider the nucleus of the elliptic surface $E(n)$, we have
that $H_2(N(n))\cong\mathbb{Z}^2$ and the corresponding
intersection matrix is the last summand in $Q_{E(n)}$, namely
$\left [ 
\begin{array}{cc}
0 & 1\\
1 & -n
\end{array}
\right ]$.
\end{subsection}

\begin{subsection}{Other building blocks.}\label{blocks}
Other ``building blocks'' are obtained as follows \cite{go95}:
\begin{itemize}
\item In $T^2\times T^2$ consider the union of the two tori 
$T^2\times p\cup
  p\times T^2$. 
Symplectically resolve the singularity and then 
blow up twice to obtain a symplectic genus $2$
  surface $F_2$ with square $0$ in $Q=T^4\#
  2\overline{\mathbb{C}\mathbb{P}^2}$ . Then $Q$ contains a symplectic 
torus $F$, disjoint from $F_2$. This $F$ is obtained from
\begin{equation*} 
F'=\{(x_1, x_2, x_3, x_4)\in \mathbb{T}^4=\mathbb{R}^4/\mathbb{Z}^4
\,|\,x_2=x_4=0\},
\end{equation*}
which is in fact Lagrangian, by perturbing the symplectic
form. Then $F$ will be disjoint from $F_2$ provided $p=(0, c)$ with $c\neq 0$.
\item Let $p$ and $q$ be distinct points in $T^2$ and consider the two 
tori $T^2\times p$ and $q\times T^2$ in $T^4=T^2\times T^2$. 
Blow up the intersection point $(q, p)$ to obtain two disjoint tori with
square $-1$, then blow up $16$ more times to reduce the square of both
to $-9$. Now take the symplectic sum of the resulting manifold $T^4\#
17\overline{\mathbb{C}\mathbb{P}^2}$ with $2$ copies of 
$\mathbb{C}\mathbb{P}^2$ along cubic curves: denote the final result 
of these operations by $S$. Then $S$ is a simply connected, symplectic
$4$-manifold, containing disjoint symplectically embedded surfaces of
genus $1$ and $2$ with trivial normal bundle. 
\item Consider a curve of degree $4$ with one transverse double point in
$\mathbb{C}\mathbb{P}^2$. The genus of such a curve is 
$2$. 
We can get a smooth surface by blowing up the the double point:
this surface represents the homology class $4h-2e$, hence it has square 
$12$. We thus need to blow up $12$ more times to get a smooth submanifold 
with genus $2$ and square $0$ in $\mathbb{C}\mathbb{P}^2\#
13\overline{\mathbb{C}\mathbb{P}^2}$.
Finally blow up three extra points, away from $F_2$, to get the manifold
$P\cong \mathbb{C}\mathbb{P}^2\#16\,\overline{\mathbb{C}\mathbb{P}^2}$.
$P$ is a symplectic simply connected manifold with Chern numbers
 $c_1^2=-7$ and $c_2=19$.
\end{itemize}  

\end{subsection}

\end{section}

\begin{section}{Construction of the examples.}

\begin{subsection}{Symplectic sphere bundles.}
\medskip

Let $(N, \beta)$ be a closed symplectic four-dimensional manifold, 
for example, one of the above building blocks, and $E\ra N$ a complex line
bundle over $N$. 
Let $\epsilon$ denote the trivial complex line bundle over $N$. 
Consider the bundle $\rho :S\ra N$ 
with fibre $S^2$ over $N$,
obtained by projectifying the complex rank two bundle $E\oplus \epsilon$.
If we denote by $E^0$ the $D^2$-bundle associated to $E$ and by 
$\overline{E^0}$ the bundle $E^0$ with opposite orientation, we can think 
of $S$ as obtained from the boundary sum $E^0\cup_{\partial E^0}
\overline{E^0}$.
Let $l$ be the tautological line bundle over $S$, as in the diagram
$$\begin{CD}
 l\subset \rho^*(E\oplus \epsilon) @>>> E\oplus \epsilon \\
 @VVV  @VVV\\
 S@>\rho>>  N 
 \end{CD}$$
If we set $c_1(l^*)=:\xi$,
there exists a ring isomorphism \cite[p. 270]{botu82}
\begin{equation*}
H^*(S;\mathbb{Z})\stackrel{\cong}{\ra}H^*(N;\mathbb{Z})[\xi]/\langle \xi^2+
\pi^*c_1(E)\xi\rangle.
\end{equation*}

Let $E^{(1)}$ denote the quotient bundle $\rho^*(E\oplus \epsilon)/l_E$.
Then the bundle $TS$ fits into the short exact sequence \cite{mu98}
\begin{equation*}
0\lra E^{(1)}\otimes l_E^*\stackrel{\alpha}{\lra} TS
\stackrel{\beta}{\lra} \rho^*TN \lra 0,
\end{equation*} 
so that the Chern classes of $S$ are given by 
\begin{eqnarray*}
c_1(S) & = & \rho^*(c_1(N)+c_1(E))+2\xi, \\
c_2(S) & = & \rho^*(c_1(N)\cup c_1(E)+ c_2(N))+2\rho^*c_1(TN)\,\xi, \\
c_3(S) & = & 2\rho^*c_2(N)\,\xi.
\end{eqnarray*}  
From this and the ring structure of $S$ we can compute the 
corresponding Chern numbers: 
\begin{eqnarray*}
c_1^3[S] & = & 6c_1^2[N]+2\langle c_1^2(E),[N]\rangle \\
c_1c_2[S] & = & 2(c_1^2[N]+c_2[N])\\
c_3[S]& = & 2c_2[N].
\end{eqnarray*}

The examples we will consider are eight-dimensional symplectic manifolds of 
the form $M=S\times F$, with $F$ a compact Riemann surface of genus $g$.
Using a product formula we can easily compute the Chern numbers of $M$:
\begin{eqnarray*}
c_4[M] & = & 2(1-g)c_3[S]=4(1-g)c_2[N]\\
c_1c_3[M] & = & 2(1-g)(c_1c_2[S]+c_3[S])
  =4(1-g)(c_1^2[N]+2c_2[N])\\
c_2^2[M] & = & 4(1-g)c_1c_2[S]=8(1-g)(c_1^2[N]+c_2[N])\\
c_1^2c_2[M] & = & 2(1-g)(c_1^3[S]+2c_1c_2[S])\\
 & = & 4(1-g)(5c_1^2[N]+2c_2[N]+\langle c_1^2(E), [N]\rangle)\\
c_1^4[M] & = & 8(1-g)c_1^3[S]=16(1-g)(3c_1^2[N]+\langle c_1^2(E), [N]\rangle).
\end{eqnarray*}
\smallskip

The next step will be to consider symplectic submanifolds of $M$ which are 
of the form $B\times \{\textrm{pt}\}$, with $B$ a symplectic submanifold 
of $S$ and $\{\textrm{pt}\}\in F$. For such submanifolds, the normal bundle 
in $M$ 
coincides with the Whitney sum of the normal bundle in $S$ and a copy of 
the trivial line bundle, which we denote again by $\epsilon$. 
This implies in particular an equivalence of 
Chern classes
\begin{equation*}
c(\nu_M B)=c(\nu_S B\oplus \epsilon)=c(\nu_S B).
\end{equation*} 

We consider for example sections $N_+$ and $N_-$ of $S$, corresponding 
to the embeddings of $N$ in $S=E^0\cup _{\partial E^0} \overline{E^0}$ as the 
zero section of $E$ and $\overline{E}$, respectively. 
In this case, the characteristic numbers which appear in the blow-up 
formulae \eqref{parpoint}, \eqref{parcurve}, \eqref{parfourdim} are given by 
\begin{eqnarray*}
c_1^2[N_{\pm}] & = & c_1^2[N]\\
c_2[N_{\pm}] & = & c_2[N]\\
\langle c_1^2(\nu_M N_{\pm}), [N_{\pm}]\rangle & = & \langle c_1^2(E), 
[N]\rangle \\
\langle c_1(\nu_M N_{\pm})c_1(N_{\pm}), [N_{\pm}]\rangle & = & \pm \langle 
c_1(N)c_1(E), [E]\rangle .
\end{eqnarray*}

Let $s$ be any such section and assume 
that $F$ is a symplectically embedded curve in $N$: then it lifts
along $s$ to a symplectically embedded curve in
$S$.
Moreover, the square of $F$ will change by an amount equal to the product
$\langle c_1(E),[F]\rangle$.
By stretching the terminology, we call here square of $F$ also the number 
resulting from evaluating the 
first Chern class of the normal bundle of $F$ (or rather, $s(F)$) in $S$ 
on its fundamental homology class. 
More precisely we have:

\begin{lemma}\label{square}
In the situation described above, the square of the lift of an embedded curve
$F$ is given by 
\begin{equation*}
\langle c_1(\nu(s(F),S)),[s(F)]\rangle=\langle c_1(\nu(F,N)),[F]\rangle
+\langle c_1(E),[F]_N\rangle,
\end{equation*}
where $\nu(\cdot ,\cdot)$ denotes the normal bundle of an embedding
and $s$ is the section under consideration.
\end{lemma}

\begin{proof}
We refer to the following commutative diagram for the notation:
$$\begin{CD}
 s(F) @>j>> S \\
 @V\cong VV  @VV\rho V\\
 F @>i>>  N 
 \end{CD}$$
We have an isomorphism of vector bundles:
\begin{align*}
\nu(s(F),S) & \cong \nu(s(F),s(N))\oplus \nu(s(N),S)|_{s(F)}\\
 & \cong \rho^* \nu(F,N)\oplus \rho^*E|_{s(F)},
\end{align*}
which implies a corresponding equivalence on the level of cohomology
classes, namely:
\begin{equation*}
c_1(\nu(s(F),S))=\rho^*c_1(\nu(F,N))+\rho^*i^*c_1(E).
\end{equation*}
We now evaluate on $[s(F)]$ and get
\begin{eqnarray*}
\langle c_1(\nu(s(F),S)),[s(F)]\rangle & = & 
 \langle \rho^*c_1(\nu(F,N)),[s(F)]\rangle+ 
    \langle \rho^*i^*c_1(E),[s(F)]\rangle\\
& = & \langle \rho^*c_1(\nu(F,N)),s_*[F]\rangle+ 
    \langle \rho^*i^*c_1(E),s_*[F]\rangle\\
& = & \langle s^*\rho^*c_1(\nu(F,N)),[F]\rangle+ 
    \langle s^*\rho^*i^*c_1(E),[F]\rangle\\
& = & \langle c_1(\nu(F,N)),[F]\rangle+ 
    \langle i^*c_1(E),[F]\rangle\\
& = & \langle c_1(\nu(F,N)),[F]\rangle+ 
    \langle c_1(E),[F]_N\rangle.
\end{eqnarray*}

Notice that, in particular, if $c_1(E)\cap [F]=0$, then the lifted curve 
has the same square in $S$ as the original one in $N$, i.e., 
\begin{equation*}
\langle c_1(\nu(s(F),S)),[s(F)]\rangle=\langle c_1(\nu(F,N)),[F]\rangle.
\end{equation*}
\end{proof}

The class $\xi\in H^2(S; \mathbb{Z})$ restricts to the standard 
K\"ahler form on each fibre of $S$. Thurston's theorem on symplectic 
fibrations \cite{mdsa98} implies that $S$ admits a symplectic form
$\omega_K=K\rho^*\beta+\eta$, which represents the class
$K\rho^*[\beta]+\xi$. 
The next lemma shows that we may in fact assume the form 
$\omega_K$ to be integral. 

\begin{lemma}\label{integral}
Given a symplectic manifold $(N, \beta)$, there exists an integral
symplectic form $\bar{\beta}$ on $N$, inducing the same Chern classes
as $\beta$. 
\end{lemma}

\begin{proof}
First we approximate $\beta $ by a closed rational form $\beta'$.
In order to do this, choose a basis $u_1, \ldots, u_m$ for $H^2(N; \mathbb{Z})$
and $2$-forms $\alpha_j\in \Omega^2(N)$ representing the element of the basis,
that is, $[\alpha_j]=u_j$. Then there exist coefficients $\lambda_j\in 
\mathbb{R}$ such that $[\beta]=\sum_{j=1}^{m}\lambda_j\alpha_j\in 
H^2(N; \mathbb{R})$.
Now consider the form 
\begin{equation*}
\beta'=\beta+\sum_{j=1}^{m}(r_j-\lambda_j)\alpha_j,\ r_j\in \mathbb{Q}.
\end{equation*}
By choosing the $r_j$'s to be rational we obtain a rational form.
In fact,
\begin{equation*}
[\beta']=[\sum_{j=1}^m r_j\alpha_j]=\sum_{j=1}^m r_ju_j\in H^2(N; \mathbb{Q}).
\end{equation*}
The differences $(r_j-\lambda_j)$ can be made arbitrarily small and for a 
sufficiently small perturbation the form $\beta'$ is still symplectic.
Moreover, since we can obviously linearly interpolate between $\beta$
and $\beta'$, the two forms induce the same Chern classes and hence 
the same Chern numbers.
Now choose a positive integer $x\in
\mathbb{Z}_{>0}$ such that $x[\beta']\in H^2(N; \mathbb{Z})$ and set
$\bar{\beta}:=x\beta'$. By construction, the form $\bar{\beta}$
is symplectic and represents an integral cohomology class. 
It is homotopic to $\beta'$, hence also to the original form $\beta$,
so it induces the same Chern classes.
\end{proof}

We may thus 
replace $\omega_K$ by ${\bar{\omega}}_K:=K'\rho^*\bar{\beta}+\eta$, 
where $K'$ is an integer larger than $K$, $\bar{\beta}$ is an integral 
symplectic form on $N$, satisfying the condition in the lemma, and $\eta$
has been chosen among the representants of $c_1(l_E^*)$. Replacing
$\omega_K$ with the new symplectic form does not affect the Chern classes,
and $[{\bar{\omega}}_K]$ is integral by construction. 

Also by construction, the forms $\omega_K$ and ${\bar{\omega}}_K$
are homotopic, so they tame the same almost complex structure.
This, together with the following lemma,
implies that symplectically embedded curves 
($2$-dimensional submanifolds) also remain symplectically embedded with 
respect to the new integral symplectic form. 

\begin{lemma}\label{pseudo4}
A smooth $2$-dimensional submanifold $F$ of a symplectic
$4$-manifold $(N, \omega)$ is symplectically
embedded if and only if it is a $J$-holomorphic curve with respect to some
tame almost complex structure $J$ on $(N, \omega)$.
\end{lemma}

\begin{proof}
Suppose first that the inclusion $i:(F, j)\ra (N, J)$ is a J-holomorphic map,
 that is, $di\circ j=J\circ di$. Let $v\in TF$.
If $i^*\omega (v, w)=0$ for all $w\in TF$, in particular 
\begin{eqnarray*}
0 & = & i^*\omega(v, jv)=\omega(di(v), di\circ j(v))\\
& = & \omega(di(v), J\circ di(v)),
\end{eqnarray*}
but the latter is strictly positive by the taming condition unless $v=0$.
Hence $i^*\omega$ is nondegenerate on $TF$, i.e., it is a symplectic form.

Conversely, suppose $F$ is a symplectic submanifold. Then by 
\cite[Lemma 3.3]{pa04}
 there exists a tame almost complex structure $J_N$ on $N$ 
such that its restriction to $TF$ is again an almost complex structure:
in fact, ${J_N}_{|TF}$ must be homotopic to $j$, since there 
exists only one homotopy 
class of almost complex structures on every orientable surface 
($\textrm{SO}(2)=\textrm{U}(1)$). 
It is then easy to perturb $J_N$ so that in fact ${J_N}_{|TF}=j$.
\end{proof}

With an integral symplectic form at our disposal, we may apply Donaldson's
existence theorem \cite{do96} 
and obtain a whole family of symplectic submanifolds 
$\{X_{\lambda}\}$ of $S$ which, for sufficiently large 
$\lambda\in \mathbb{Z}$, 
realise the Poincar\'e dual of $\lambda[\omega_K]$, i.e.,
$PD_S[X_{\lambda}]=\lambda(\xi+K[\pi^*\beta])$ in $H^2(S;\mathbb{Z})$. 

Let $i$ denote the inclusion $X_{\lambda}$ in $M$. We have the relation 
$$i^*c(S)=c(X_{\lambda})\cup\, c(\nu_S X_{\lambda})$$
and since $c_1(\nu_S X_{\lambda})=e(\nu_M X_{\lambda})=i^*PD_M[X_{\lambda}]$,
 we can rewrite it as 
$$i^*c(M)=c(X_{\lambda})\cup\, (1+i^*\lambda(\xi+K[\pi^*\beta])).$$
From this relation, using 
\begin{equation*}
\langle i^*y, [X_{\lambda}]\rangle =\langle y\cup(\lambda\xi+
\lambda K[\pi^*\beta]), [M]\rangle\ \textrm{for\ all}\  y\in H^4(M)
\end{equation*}
 and the cohomology ring structure of $M$, we can
compute the invariants of $X_{\lambda}$ (see \cite[Appendix A]{pa04} 
for explicit computations).  
They are:
\begin{eqnarray*}
c_2[X_{\lambda}]&=& \lambda c_2[N]+\lambda (\lambda -1) \langle
c_1(N)c_1(L),[N]\rangle \nonumber\\
& & -2\lambda K (\lambda -1) \langle c_1(N)[\beta],[N]\rangle
  +\lambda^2(\lambda-1)\langle c_1^2(E),[N]\rangle \nonumber\\
& & +\lambda^2K(2-3\lambda)\langle c_1(E)[\beta],[N]\rangle\nonumber\\
& &  +\lambda^2K^2(3\lambda -2)\langle [\beta]^2,[N]\rangle \nonumber
\end{eqnarray*}
\begin{eqnarray} \label{donald}
c_1^2[X_{\lambda}]&=& \lambda c_1^2[N]+2\lambda (\lambda -1)
\langle c_1(N)c_1(L),[N]\rangle\nonumber\\
& & +4\lambda K(1-\lambda)\langle c_1(N)[\beta],[N]\rangle\nonumber\\
& & +\lambda(\lambda^2-2\lambda+1)\langle c_1^2(E),[N]\rangle \\
& &+\lambda^2K(4-3\lambda)\langle c_1(E)[\beta],[N]\rangle\nonumber\\
& &  +\lambda^2K^2(3\lambda-4)\langle [\beta]^2,[N]\rangle\nonumber\\
\langle c_1^2(\nu_M X_{\lambda}),[X_{\lambda}]\rangle & = & 
\lambda^3\langle c_1^2(E),[N]\rangle
-3\lambda^3K\langle c_1(E)[\beta],[N]\rangle \nonumber\\
& & +3\lambda^3K^2\langle [\beta]^2,[N]\rangle \nonumber\\
\langle c_1(X_{\lambda})c_1(\nu),[X_{\lambda}]\rangle & = &  
-\lambda^2\langle c_1(N)c_1(L),[N]\rangle+2\lambda^2K
\langle c_1(N)[\beta],[N]\rangle \nonumber\\
& & +\lambda^2(1-\lambda)\langle c_1^2(E),[N]\rangle
+\lambda^2K(3\lambda-2)\langle c_1(E)[\beta],[N]\rangle \nonumber\\
& & +\lambda^2K^2(2-3\lambda)\langle [\beta]^2,[N]\rangle. \nonumber
\end{eqnarray}
\end{subsection}
\end{section}

\begin{section}{The blow-up systems.}

Now we have introduced all the elements necessary for the proof of Theorem 
\ref{geopar}.
The proof itself follows that of Halic for the case of dimension $6$.
For any given $j$, that is, we will show that it is possible to construct
a symplectic manifold realising $j$ and with enough symplectic submanifolds
so that we can vary the other parameters and eventually realise all
admissible quintuples. We construct our examples distinguishing 
three main cases. 
 
\begin{subsection}{Realising sets of parameters with $j\geq 1$.}
\medskip

We start by considering examples of symplectic $8$-dimensional manifolds
for which the parameter $j$ is greater than or equal to $1$. 
In order to produce such examples, we perform the symplectic sum of
the manifolds $Q$ and $E(n)$ of paragraph \ref{blocks} along 
symplectically embedded tori with square zero. Before doing so, 
though, we blow up one
extra point in $Q$, away from the torus along which we intend to
perform the sum.

The result of these operations is the manifold 
$$X_n:=Q^*\#_F E(n),$$
which is a simply connected symplectic manifold, with Chern numbers
 $c_1^2=-3$ and $c_2=3+12n$.

Following \cite[Prop. 2.2]{sm00}, we may write
 $$c_1(X_n)=c_1(Q^*)+c_1(E(n))-2PD([F]).$$

Assume that the torus $F\subset E(n)$ is actually contained in $N(n)$, the
nucleus of the elliptic surface (cf. Definition \ref{nucleus}), and recall 
that we denote by
$\tau_{ij}$, $i=1, \ldots , n$, $j=1, \ldots , 8$, the elements of the 
basis of $H_2(E(n))$ corresponding
to the $n$ copies of the $(-E_8)$-block in the intersection matrix 
(\ref{Enint}).
The $\tau_{ij}$'s are represented by submanifolds of 
the complement of $N(n)$, hence
disjoint from $F$: for this reason they represent homology classes in $X_n$ 
(which we still denote by $\tau_{ij}$) and we may consider their
Poincar\'e duals, which will be elements of $H^2(X_n)$.

According to this interpretation of the elements $\tau_{ij}\in H^2(X_n)$
we consider the complex line bundle $L$ over $X_n$, specified 
by its first Chern class 
$$c_1(L)=\sum_{i=1}^{n-1}2PD(\tau_{i1})+\sum_{j=1}^3PD(\tau_{nj})+
PD(\tau_{n8})\in H^2(X_n; \mathbb{Z}),$$
and compute 
$$\left \{
\begin{array}{l}
\langle c_1^2(L), [X_n]\rangle =-8(n-1)-12\\
\langle c_1(L)c_1(X_n), [X_n]\rangle =0.
\end{array}
\right.$$

Now consider the manifold $S=\mathbb{P}(L\oplus \mathbb{C})$.
Let $s:X_n\ra S$ be the section which embeds $X_n$ in $S$ as
$\mathbb{P}(L\oplus \{0\})=({X_n})_-$. Denote by $E$ the exceptional 
sphere of
  the last blow-up in $Q^*$ and by $e=[E]$ its fundamental homology class. 
Similarly,
  if $E_-=s(E)$ is the lift of $E$ along the section $s$, let
  $e_-=[E_-]=[s(E)]$ be the corresponding homology class, so that we have
  the following commutative diagram:
\[
\begin{CD}
 E_- @>>> S \\
 @V\rho V\cong V  @VV\rho V\\
 E @>>>  X_n 
 \end{CD}
\]
Then
$\langle c_1(L_{|E}), e\rangle =PD^{-1}c_1(L)\cdot e=0$, so by Lemma
\ref{square} we have
\begin{equation*}
\langle c_1(\nu_M E_-), e_-\rangle =e^2+\langle c_1(L_{|E}), e\rangle =e^2.
\end{equation*}
Observe that if $F_2$ denotes the symplectically embedded surface of genus
$2$ in $Q^*$, it is disjoint from any representative of the classes
$[\tau_{ij}]$ and therefore it also lifts to ${F_2}_-\subset {X_n}_-$ with
\begin{equation*}
\langle c_1(\nu_M {F_2}_-), [{F_2}]_-\rangle =0.
 \end{equation*}
Finally recall that $S$ admits an integral symplectic form $\omega_K$ 
and hence, for
$\lambda$ large enough, symplectic submanifolds $X_{\lambda}$, realising the
Poincar\'e dual of $\lambda[\omega_K]$, whose invariants are given by the
expressions in \ref{donald}.

Take the product of $S$ with $S^2$ to obtain the simply connected
symplectic
$8$-dimensional manifold $M=S\times S^2$. The parameter $j$ (recall: 
$j$ was defined by the congruence relation
$-c_4+c_1c_3+c_2^2+4c_2c_1^2-c_1^4=720j$) in this case 
takes on precisely 
the value $n$. The other parameters are 
\begin{eqnarray*}
a & = & 48n+12\\
4m & = & -12\\
12k & = & -192n-468\\
b & = & -128n-208.
\end{eqnarray*}
 
Now blow up $M$ at $x$ points, $y$ copies of $E$, $z$ copies of $F_2$, 
$u$ copies of ${X_n}_-$ and $v$ copies of $X_{\lambda}$. Denote by 
$\widetilde{M}$
the manifold obtained after performing these blow-ups and let $(a', m', j',
k', b')$ be the parameters associated to $\widetilde{M}$. Then $j'=j$
(recall that we already remarked that $j$ is invariant under blow-up),
whereas by applying the blow-up formulae (\ref{parpoint}), 
(\ref{parcurve}) and (\ref{parfourdim})
we find the following expressions for the other parameters:
\begin{eqnarray} \label{blup}
a' & = & 48n+12+3x+4y-4z+(12n+3)u+b_1v\nonumber\\
4m' & = & -12-4y+4z+(-36n-12)u+b_2v\\
12k' & = & -192n-468-180x-108y+144z+(36n+60)u+b_3v\nonumber\\
b' & = & -128n-208-81x-48y+64z+(24n+30)u+b_4v\nonumber.
\end{eqnarray}
The coefficients coming from blowing up along the submanifold
$X_{\lambda}$, whose invariants are computed in the Appendix,
are given by:
\begin{eqnarray*}
b_1 & = & \lambda
c_2(N)+(\lambda^3-\lambda^2)\,c_1^2(E)+(3\lambda^3K^2-2\lambda^2K^2)\,
[\beta]^2\\
& & +(\lambda^2-\lambda)\,c_1(N)\cup\,c_1(E)
 +(2\lambda K-2\lambda^2K)\,c_1(N)\cup\,[\beta]\\
& & +(2\lambda^2K-3\lambda^3K)\,c_1(E)\cup\,[\beta]\\
b_2 & = & \lambda c_1^2(N)-3\lambda
c_2(N)+(\lambda^2-2\lambda^3+\lambda)\,c_1^2(E)+
(2\lambda^2-6\lambda^3K^2)\,[\beta]^2\\
 & & +(\lambda-\lambda^2)\,c_1(N)\cup\,c_1(E)
+2(\lambda^2 K-\lambda K)\,c_1(N)\cup\,[\beta]\\
 & & +(6\lambda^3 K-2\lambda^2
 K)\,c_1(E)\cup\,[\beta]\\
b_3 &  = & -13\lambda c_1^2(N)-\lambda
c_2(N)+(9\lambda^2-2\lambda^3-13\lambda)\,c_1^2(E)\\
 & & +(18\lambda^2K^2-6\lambda^3 K^2)\,[\beta]^2+
  (27\lambda-9\lambda^2)\,c_1(N)\cup\,c_1(E)\\
 & & +(18\lambda^2 K-54\lambda K)\,c_1(N)\cup\,[\beta]
  +(6\lambda^3
 K-18\lambda^2 K)\,c_1(E)\cup\,[\beta] \\
b_4 & = & -6\lambda c_1^2(N)+(4\lambda^2-6\lambda-\lambda^3)\,c_1^2(E)+
(8\lambda^2
K^2-3\lambda^3 K^2)\,[\beta]^2+\\
 & & (12\lambda-4\lambda^2)\,c_1(N)\cup\,c_1(E)+(8\lambda^2 K-24\lambda
 K)\,c_1(N)\cup\,[\beta]\\
 & & +(3\lambda^3
 K-8\lambda^2 K)\,c_1(E)\cup\,[\beta],
\end{eqnarray*}
where we have suppressed evaluation on the fundamental class of $[N]$ from 
the notation.

We regard (\ref{blup}) as a linear system in the variables $x$, $y$, $z$,
$u$, $v$. If we can prove that for arbitrary parameters $a'$, $m'$, $k'$,
$b'$, satisfying the additional condition 
$a'+m'\equiv 0\,(\textrm{mod}\ 3)$, 
the system admits a quintuple of positive, integer solutions, then we
will have shown that we can realise all such parameters, precisely by 
performing on $M$ the
sequence of blow-ups corresponding to the solutions $(x, y, z, u, v)$. 

The solutions of system (\ref{blup}) are:
\begin{eqnarray*}
x & = & [(32n+13)\lambda^3 K^2 [\beta]^2+r_1(\lambda, K)]v+
\bigl(8n+\frac{10}{3}\bigr)a'+\\
 & & \bigl(32n+\frac{40}{3}\bigr)m'+(-128n-48)k'+(24n+9)b'+\\
 & & 640n^2+224n\\
y & = & \bigl[\bigl(12n+\frac{9}{2}\bigr)\lambda^3 K^2 [\beta]^2
+r_2(\lambda, K)\bigr]v
      +(3n+3)a'+\\
 & & (12n+8)m'+(-48n-21)k'+(9n+4)b'+\\
 & & 240n^2+32n+1\\
z & = & [(48n+18)\lambda^3 K^2 [\beta]^2+r_3(\lambda, K)]v+(12n+6)a'+\\
 & & (48n+21)m'+(-192n-69)k'+(36n+13)b'+\\
 & & 960n^2+272n+4
\end{eqnarray*}
\begin{eqnarray*}
u & = & [4\lambda^3 K^2 [\beta]^2+r_4(\lambda, K)]v+a'+4m'-16k'+\\
 & & 3b'+80n.
\end{eqnarray*}

First of all notice that these solutions are, indeed, integer, 
because of the additional condition $a'+m'\equiv
  0\ (\textrm{mod}\ 3)$.
Moreover, one can observe that the $r_i(\lambda, K)$ are polynomials of 
degree at most $2$ in
  $\lambda$ and $2$ in $K$, with coefficients depending on 
 $c_1^2(N), c_2(N), c_1^2(E), [\beta]^2, c_1(N)\cup\,c_1(E), 
 c_1(N)\cup\,[\beta], c_1(E)\cup\,[\beta]$, 
 evaluated on $[N]$; 
recall also that $\langle [\beta]^2, [N]\rangle $ is
  strictly positive, because $\beta^2=\beta\wedge\beta$ is a volume
form on $N$:
by the previous two remarks we conclude that, by choosing $\lambda$
  large enough, we may ensure the positivity of the $v$-coefficients and
  consequently, again by a choice of $v$ sufficiently large, positivity of
  all variables. 
\end{subsection}

\begin{subsection}{The case $j=0$.}

In order to show that all quintuples of parameters with $j=0$ admit a
symplectic realisation, we start again by constructing a
manifold with $j=0$. We consider the $4$-manifold $Q^*$ and the complex
line bundle $L$ defined by $c_1(L)=-2e_3$, with $e_3$ denoting the
exceptional divisor of the last blow-up.
Since
$$c_1(Q^*)=c_1(\mathbb{T}^4)-\sum_{i=1}^3 e_i,$$
we see that 
$$\left \{
\begin{array}{l}
\langle c_1^2(L), [Q^*]\rangle =-4\\
\langle c_1(L)c_1(X_n), [Q^*]\rangle =-2.
\end{array}
\right.$$ 

We proceed to construct $S=\mathbb{P}(L\oplus\mathbb{C})$ and 
$M=S\times S^2$ as in the previous section.
Then $M$ realises $j=0$, as required, and the other parameters are
\begin{eqnarray*}
a & = & 12\\
4m & = & -24\\
12k & = & -488\\
b & = & -218.
\end{eqnarray*} 
 
We let $Q^*_-$ be the image of the embedding of $Q^*$ in $M$ as 
$\mathbb{P}(L\oplus
\{0\})$. Then $Q^*_-$ contains a sphere $E$ with square $-1$ (the
exceptional sphere of either the first or the second blow-up) and a genus
$2$ surface $F_2$. These curves intersect in $S$, but since 
the cup product of their Poincar\'e duals and $c_1(L)$ vanishes,
by Lemma \ref{square} they provide disjoint
submanifolds $E\times \{\rm{pt}\}$ and $F_2\times \{\rm{pt}\}$ of 
$M=S\times S^2$
  with the same genus and square. 

Together with the submanifold $Q^*_-$, we consider as in the previous
cases submanifolds $X_{\lambda}$, realising the Poincar\'e duals of
multiples $\lambda\,\omega_K$ of some integral symplectic form $\omega_K$ on
 $M$.

We are now able to write down the blow-up system for $M$, where we
blow up at $x$ points, $y$ copies of $E$, $z$ copies of $F_2$, 
$u$ copies of ${Q^*}_-$ and $v$ copies of $X_{\lambda}$.
Then the parameters of our $8$-dimensional manifold transform according
to the following expressions
\begin{eqnarray} \label{blup0}
a' & = & a+3x+4y-4z+3u+b_1v\nonumber\\
4m' & = & 4m-4y+4z-12u+b_2v\\
12k' & = & 12k-180x-108y+144z+24u+b_3v\nonumber\\
b' & = & b-81x-48y+64z+14u+b_4v,\nonumber
\end{eqnarray}
where the $b_i$'s are once again the coefficients corresponding to 
blow-up
along submanifolds belonging to the family $\{X_{\lambda}\}$.

The solutions in this case are 
\begin{eqnarray*}
x & = & [(13\lambda^3 K^2 [\beta]^2+r_1(\lambda, K)]v+\frac{10}{3}a'+
  \frac{40}{3}m'-48k'+9b'\\
y & = & \bigl[\frac{17}{2}\lambda^3 K^2 [\beta]^2+r_2(\lambda, K)\bigr]
v+4a'+12m'-37k'+7b'+1\\
z & = & [22\lambda^3 K^2 [\beta]^2+r_3(\lambda, K)]v+7a'+
  25m'-85k'+16b'+4\\
u & = & [4\lambda^3 K^2 [\beta]^2+r_4(\lambda, K)]v+a'+4m'-16k'+3b'.
\end{eqnarray*}
Once again, we may observe that all solutions are integer and that by a 
suitable choice of $\lambda$ and $v$ we may assume them to be positive, 
as well.
Hence all admissible quintuples with $j=0$ admit a symplectic 
realisation.
\end{subsection}

\begin{subsection}{Negative values of $j$.}

We are left with only the case $j<0$ to take care of. For this we 
construct a $6$-dimensional manifold $S$
as in the cases of positive values of the parameter $j$ and then 
define
$M$ to be the product of $S$ with a compact Riemann surface of genus
two. Notice that in this case the realisation will not
be simply connected.

The only difference in the blow-up system occurs in the parameters
corresponding to the manifold $\Sigma$ which is blown up: these have in
fact opposite sign.
Therefore the blow-up system has the form 
\begin{eqnarray*}
a' & = & -48n-12+3x+4y-4z+(12n+3)u+b_1v\\
4m' & = & 12-4y+4z+(-36n-12)u+b_2v\\
12k' & = & 192n+468-180x-108y+144z+(36n+60)u+b_3v\\
b' & = & 128n+208-81x-48y+64z+(24n+30)u+b_4v
\end{eqnarray*}
and the solutions are given by
\begin{eqnarray*}
x & = & [(32n+13)\lambda^3 K^2 [\beta]^2+r_1(\lambda, K)]v+
(8n+\frac{10}{3})a'+\\
 & & (32n+\frac{40}{3})m'+(-128n-48)k'+(24n+9)b'+\\
 & & -640n^2-224n\\
y & = & [(12n+\frac{9}{2})\lambda^3 K^2 [\beta]^2+r_2(\lambda, K)]v+
(3n+3)a'+\\
 & & (12n+8)m'+(-48n-21)k'+(9n+4)b'+\\
 & & -240n^2-32n-1\\
z & = & [(48n+18)\lambda^3 K^2 [\beta]^2+r_3(\lambda, K)]v+(12n+6)a'+\\
 & & (48n+21)m'+(-192n-69)k'+(36n+13)b'+\\
 & & -960n^2-272n-4\\
u & = & [4\lambda^3 K^2 [\beta]^2+r_4(\lambda, K)]v+a'+4m'-16k'+\\
 & & 3b'-80n.
\end{eqnarray*}
The same considerations as to positivity and integrality apply as in the case
of positive $j$.
\end{subsection}
\end{section}

\begin{section}{Some final remarks.}
\begin{subsection}{K\"ahler manifolds.}

Symplectic manifolds occupy the central position in the sequence 
of inclusions
\begin{equation*}
\textrm{K\"ahler}\subsetneq\textrm{symplectic}\subsetneq
\textrm{almost complex}.
\end{equation*} 
These inclusions have long been known to be proper. 
One is interested in finding out which properties distinguish 
symplectic manifolds from the manifolds in the other two classes and
which ones do not. 

We have shown in this paper that in dimension $8$,
the geography of symplectic manifolds
coincides with that of almost complex manifolds.
It is then natural to ask whether this is true also for
the geography of K\"ahler manifolds.
Our work unfortunately does not provide a positive answer to 
this question. 
The examples we construct, in fact, might already cease to be K\"ahler 
at the point where we take symplectic sums of
$4$-dimensional building blocks. 
\end{subsection}

\begin{subsection}{Chern numbers and topology.}
The Chern numbers of a symplectic manifold
are invariants of the symplectic form. 
When considering invariants of some nature, it is always interesting and 
natural to ask to what extent these invariants classify.
The Chern numbers certainly fail to classify symplectic structures.
Deformation equivalent symplectic forms have the same Chern numbers and so
do isomorphic symplectic forms.
There are even examples of symplectic forms which are not related by any 
sequence of
isomorphisms and deformation equivalences and which are not distinguished 
by the Chern numbers.
In other words, the extent to which Chern numbers fail to classify symplectic
structures is considerable. Therefore one may wonder
whether they might not be topological invariants.

In dimension $4$, the Chern numbers $c_1^2$ and $c_2$ are indeed 
topological invariants.
In dimension $6$, LeBrun has shown that Chern numbers are not 
topological invariants of complex manifolds, but what happens if
we introduce a symplectic form is not known. 
A better understanding of the geography of symplectic manifolds
may be useful in order to answer this question by 
comparing the Chern numbers with the topology.
 
The symplectic constructions on which Halic's results and our main theorem
rely allow a good control of the cohomological data.
In dimension, in particular, we have the smooth classification 
result of Wall \cite{wa66} based on those data.
So one might hope to be able to detect a smooth manifold realising two 
different combinations of Chern numbers, that is, admitting two distinct 
symplectic structures, distinguished by the Chern numbers.
\end{subsection}

\begin{subsection}{Geography with fundamental group.}\label{geoG}
We would like to conclude by briefly addressing the question of geography 
with fundamental group, that is: to which extent is it possible to prescribe 
Chern numbers and fundamental group of a symplectic manifold at the same time.

Observe that Halic' result yields in dimension $6$ simply connected 
realisations for all admissible triples.
In dimension $4$, on the other hand, if $(p, q)$ is an admissible pair 
with $p+q<0$, there exists no simply connected symplectic manifold 
with $(c_1^2, c_2)=(p, q)$ (\cite[Prop. 4.20]{pa04}).

If $G$ is any finitely presentable group, Gompf has shown that there exists 
a closed symplectic $4$-manifold with fundamental group $G$.
This manifold, moreover, may be assumed to satisfy certain additional 
properties.

\begin{thm}[Gompf]\label{gompf}
Let $G$ be any finitely presentable group. Then there is a closed 
symplectic manifold $M_G$ with $\pi_1(M_G)\cong G$. Furthermore we may 
assume:
\begin{itemize}
\item [{\rm(i)}] $c_1^2[M_G]=0$, $c_2[M_G]=12r>0$;
\item [{\rm(ii)}] $M_G$ contains a symplectic torus $T$ with square $0$
and inclusion $i:T\ra M_G$ inducing the trivial map on $\pi_1$.
\end{itemize}
\end{thm}
The proof can be found in \cite{go95} and relies in fact 
on the symplectic connected sum construction.

\begin{rem}
The existence of a four-dimensional symplectic manifold $M_G$ with
$\pi_1(M_G)\cong G$ for every finitely presentable group $G$ is another 
feature that distinguishes symplectic from K\"ahler manifolds.
The abelianisation of the fundamental group of a K\"ahler manifold, 
in fact, necessarily has even rank, since $(\pi_1)_{\textrm{ab}}=
\pi_1/[\pi_1, \pi_1]\cong H_1$.
\end{rem}
   
Theorem \ref{gompf} can be applied to improve partially Halic' result 
in dimension $6$ and show that some admissible triples may be realised
by a closed connected symplectic manifold $M$ having a given 
finitely presentable fundamental group $G$.

\begin{prop}
For every admissible triple $(2a, 24b, 2c)$ with $b\leq -2$ and every 
finitely presentable group $G$ there exists a closed symplectic $6$-manifold
$M$ such that 
\begin{eqnarray*}
c_1^3[M] & = & 2a\\
c_1c_2[M] & = & 24b\\
c_3[M] & = & 2c\\
\pi_1(M) & \cong & G.
\end{eqnarray*} 
\end{prop}

\begin{proof}
Let $S$ be as in Section \ref{blocks} and denote by $X$ the symplectic 
connected sum of $S$ with $E_n^*$, the blow-up at one point of the 
elliptic surface $E_n$.
In order to realise all admissible triples with $b\leq -2$, Halic considers 
the manifold:
\begin{equation*}
M'=X\times F_2\#_{F_1\times F_2}S\times F_1.
\end{equation*}
For a given finitely presented group $G$, let $M$ be the manifold 
obtained by taking the connected sum of $M'$ with the product 
$M_G\times F_1$, where $M_G$ is as in Theorem \ref{gompf} and
$F_1$ is a surface of genus $1$. In other words we have
\begin{equation*}
M=M'\#_{F_1\times F_1} M_G\times F_1,
\end{equation*}
where the sum is performed on the left-hand side along 
$F_1\times F_1\subset S\times F_1$ (this is possible because 
$F_1$ and $F_2$ are disjoint in $S$) and $F_1=T\subset M_G$ is
also as in (ii) of Theorem \ref{gompf}.
The Chern numbers of $M$ are the same as those of $M'$ by 
\cite[Prop. 1.3]{ha99}. Moreover, one can show that 
the fundamental group $\pi_1(M)$ is isomorphic to 
$G$. To see that, let
\begin{equation*}
U=M'-F_1\times F_1\ \textrm{and}\ V=M_G\times F_1-F_1\times F_1.
\end{equation*}
Then $\{U, V\}$ is an open covering of $M$ and 
$U\cap V=F_1\times F_1\times A$. 
Observe that we may rewrite $U$ as 
\begin{equation*}
X\times F_2\#_{F_1\times F_2} (S-F_1)\times F_1.
\end{equation*}
Since the complement of $F_1$ in $S$ is simply connected, the argument 
in \cite[p.~379]{ha99} still goes through and shows that $U$ is 
simply connected.
By Seifert-van Kampen, the fundamental group of $M$ is given by
\begin{equation}\label{SvK}
\pi_1(M)=\pi_1(V)/\langle\pi_1(U\cap V)\rangle =\pi_1(M_G-F_1)/
\langle\pi_1(F_1)\times \pi_1(A)\rangle.
\end{equation}
The epimorphism
\begin{equation*}
i_*:\pi_1(M_G-F_1)\ra \pi_1(M_G),
\end{equation*}
is surjective because of the codimension of the embedding 
$F_1\subset M_G$. Moreover, the kernel of $i_*$ is generated by a meridian 
of $F_1$ and can be identified with $\pi_1(A)$, so we have 
\begin{equation*}
\pi_1(M_G-F_1)/\langle\pi_1(F_1\times A)\rangle\cong
\pi_1(M_G)/\langle\pi_1(F_1)\rangle\cong G,
\end{equation*}
which implies, together with (\ref{SvK}), that $\pi_1(M)\cong G$.
By blowing up $M$ we get symplectic realisations with fundamental 
group $G$ for all admissible triples with $b\leq -2$.
\end{proof}

\end{subsection}
\end{section}

\end{document}